\documentclass[a4paper,11pt,twoside]{article}
\usepackage{pst-fill,pst-grad}
\usepackage{textcomp}
\usepackage[english]{babel}
\usepackage{graphicx}
\usepackage{amsmath}
\usepackage{float}
\usepackage[matrix,arrow,curve]{xy}
\usepackage{pstricks}
\usepackage{amsmath,amsfonts,verbatim,afterpage,theorem,euscript,mathrsfs,amssymb}
\usepackage{amsfonts}
\usepackage{amssymb}
\usepackage{array}
\usepackage{dsfont}
\usepackage{hyperref}
\newtheorem{Theoreme}{Theorem}

\newtheorem{Lemme}{Lemma}[section]
\newtheorem{Corollaire}{Corollary}[section]
\newtheorem{Remarque}{Remark}[section]
\title{\bf Functional Inequalities in Stratified Lie groups with Sobolev, Besov, Lorentz and Morrey spaces}
\author{Diego Chamorro\footnote{Laboratoire de Math\'ematiques et Mod\'elisation d'Evry (LaMME) - UMR 8071. Universit\'e d'Evry Val d'Essonne, 23 Boulevard de France, 91037 Evry Cedex, France. email: \textit{diego.chamorro@univ-evry.fr}}, Anca-Nicoleta Marcoci\footnote{Department of Mathematics and Computer Science. Technical University of Civil Engineering, Bucharest, Bld. Lacul Tei, no. 124, sector 2. Romania. email: \textit{anca.marcoci@utcb.ro}}, Liviu-Gabriel Marcoci\footnote{Department of Mathematics and Computer Science. Technical University of Civil Engineering, Bucharest, Bld. Lacul Tei, no. 124, sector 2. Romania.  email: \textit{lmarcoci@instal.utcb.ro}}.}
\DeclareMathOperator*{\esssup}{ess\,sup}
\begin{document}
\maketitle
\begin{scriptsize}
\abstract{ When $p>1$, using as base space classical Lorentz spaces associated to a weight from the Ari\~no-Muckenhoupt class $B_p$, we will study Gagliardo-Nirenberg inequalities. As a by-product we will also consider Morrey-Sobolev inequalities. These arguments can be generalized to many different frameworks, in particular the proofs are given in the setting of stratified Lie groups.}\\[3mm]
\end{scriptsize}
\section{Introduction and presentation of the results}
The aim of this article is to provide, in the setting of stratified Lie groups, a general proof for a particular type of improved Sobolev inequalities of the following general form
\begin{equation}\label{Inequality_1}
\|f\|_{\dot{W}^{s_1,q}}\leq C \|f\|_{\dot{W}^{s,p}}^\theta\|f\|_{\dot{B}^{-\beta,\infty}_{\infty}}^{1-\theta},
\end{equation}
where $f\in \dot{W}^{s,p} \cap \dot{B}^{-\beta,\infty}_{\infty}$. Here we write $\dot{W}^{s,p}$ for homogeneous Sobolev spaces and $\dot{B}^{-\beta,\infty}_{\infty}$ for homogeneous Besov spaces (see Section \ref{Secc_Functional_Spaces} below for precise definitions). The parameters $s, s_1, p, q$ and $\beta$ defining Sobolev and Besov spaces in the previous inequality are related by the conditions $1<p<q<+\infty$, $\theta=p/q$, $s_1=\theta s -(1-\theta)\beta$ and $-\beta<s_1<s$, but they do not depend on the dimension and in this sense these inequalities are more general than classical Sobolev inequalities; of course the inequalities above are sharper than classical ones. Historically, the first proof in the Euclidean setting of these inequalities is due to P. G\'erard, F. Oru and Y. Meyer \cite{GMO} and is based on a Littlewood-Paley decomposition and interpolation results applied to dyadic blocks. Another proof of these inequalities using maximal function and Hedberg's inequality is given in \cite{Chamorro1}.\\ 

Let us mention now that in the case of Lorentz spaces $L^{p,q}$, H. Bahouri and A. Cohen \cite{Bahouri} proved the inequality
\begin{equation}\label{Inequality_6}
\|f\|_{L^{p,q}}\leq C \|f\|_{\dot{B}^{s,q}_q}^{q/p}\|f\|_{\dot{B}^{s-n/q,\infty}_q}^{1-q/p} \quad \mbox{ with } \frac{1}{p}=\frac{1}{q}-\frac{s}{n}.
\end{equation}
Remark that in this estimate, the index $q$ defining the Lorentz and the Besov spaces is related to the parameters $p$, $s$ and the dimension $n$. This inequality was generalized in the Euclidean setting by D. Chamorro \& P-G. Lemari\'e-Rieusset \cite{Chamorro3} to other values of the parameter $q$ using interpolation techniques and pointwise estimates. In a recent article, V.I. Kolyada and F.J. P\'erez L\'azaro \cite{Kolyada} gave an interesting proof for inequalities of type (\ref{Inequality_1}) and (\ref{Inequality_6}) based on the use of rearrangement inequalities and the properties of the Gauss-Weierstrass kernel.\\

Motivated by the use of Lorentz spaces in these previous works, in our first theorem we will provide a generalization (in stratified Lie groups) of improved Sobolev inequalities of type (\ref{Inequality_1}) by considering weighted Lorentz-based Sobolev spaces defined as the set of measurable functions $f:\mathbb{G} \longrightarrow \mathbb{R}$ such that the quantity $\|(-\Delta)^{\frac{s}{2}}f\|_{\Lambda^p(w)}$ is bounded where for $s>0$ the fractional power of the Laplacian is defined in Section \ref{Secc_Notation} below, and for $1<p<+\infty$ the space $\Lambda^p(w)$ denotes the classical Lorentz space of functions introduced in \cite{Lorentz} and \cite{Lorentz1} defined as
$$\Lambda^p(w)=\left\{f:\|f\|_{\Lambda^p(w)}=\left(\int_0^{+\infty} f^*(t)^pw(t)dt\right)^{1/p}<+\infty\right\},$$
where $w$ is a weight in $\mathbb{R}_+$ and $f^*$ denotes the non-increasing rearrangement of $f$ (see \cite{BS} for standard notations). Many of the properties of these spaces depend on the weight $w$: in particular, if $w=1$ we have $\Lambda^p(w)=L^p$ and if $w(t)=t^{p/q-1}$, with $1\leq q\leq+\infty$, we obtain $\Lambda^p(w)=L^{q,p}$, where $L^{q,p}$ are the usual Lorentz spaces. In this work we will consider the weighted Lorentz space $\Lambda^p(w)$ such that the weight $w$ satisfies the $B_p$ condition, the reason for this is given by the fact that M. A. Ari\~no and B. Muckenhoupt showed in \cite{Arino} that this $B_p$ condition characterizes the boundedness of the Hardy-Littlewood maximal operator on $\Lambda^p(w)$ and this particular property will be intensively used in our proofs.  See Section \ref{Secc_Functional_Spaces} for definitions and \cite{Arino}, \cite{CarroSoria}, \cite{Carro} and \cite{Soria} for more details and properties concerning these functional spaces.\\

In this direction of generalization, standard Sobolev inequalities have been studied by A. Cianchi \cite{Cianchi} in the context of Orlicz-Sobolev spaces, but improved inequalities of the general type (\ref{Inequality_7}) presented in Theorem \ref{Theorem2} below are, to the best of our knowledge, new.
\begin{Theoreme}\label{Theorem2} Let $\mathbb{G}$ a stratified Lie group. Let $s>0$, $w\in B_p$ be a weight and let $f:\mathbb{G}\longrightarrow\mathbb{R}$ be a function such that $(-\Delta)^{\frac{s}{2}}f\in \Lambda^{p}(w)(\mathbb{G})$ and $f\in \dot{B}^{-\beta, \infty}_\infty(\mathbb{G})$. Then we have the following version of improved Sobolev inequalities:
\begin{equation}\label{Inequality_7}
\|(-\Delta)^{\frac{s_{1}}{2}}f\|_{\Lambda^{q}(w)}\leq C \|(-\Delta)^{\frac{s}{2}}f\|_{\Lambda^{p}(w)}^\theta \|f\|^{1-\theta}_{\dot{B}^{-\beta, \infty}_\infty},
\end{equation}
where $1<p<q<+\infty$, $\theta=p/q$, $s_1=\theta s -(1-\theta)\beta$ and $-\beta<s_1<s$.
\end{Theoreme}
The choice of the weights in the $B_p$ class is given by two important facts. First, these weights allow us to consider general functional spaces, in particular we can easily recover standard Lorentz spaces. Second, these weights ensure that maximal function is bounded in the spaces $\Lambda^{p}(w)$ for $1<p<+\infty$, and this feature is crucial as the proof of Theorem \ref{Theorem2} requires this property. Note in particular that inequality (\ref{Inequality_7}) is different from inequality (\ref{Inequality_6}) since Lorentz-Sobolev spaces are not included in the scale of Besov spaces.\\

Since our proof of Theorem \ref{Theorem2} relies essentially on a pointwise inequality and on the boundedness of the Hardy-Littlewood maximal operator, it is possible to give a related result replacing classical Lorentz spaces by Morrey spaces $\mathcal{M}^{p,a}$ which are a useful generalization of Lebesgue spaces. Classical Hardy-Littlewood-Sobolev inequalities were studied in this functional framework by D. Adams \cite{Adams} and by F. Chiarenza \& M. Frasca \cite{Chiarenza} and our next theorem is an improvement of these inequalities.
\begin{Theoreme}\label{Theorem3} Let $\mathbb{G}$ be a stratified Lie group. Let $s>0$, $1<p<+\infty$ and $0\leq a < n$ and let $f$ be a function such that $(-\Delta)^{\frac{s}{2}} f\in \mathcal{M}^{p,a}(\mathbb{G})$ and  $f\in\dot{B}^{-\beta, \infty}_\infty(\mathbb{G})$. Then we have
\begin{equation}\label{Inequality_8}
\|(-\Delta)^{\frac{s_{1}}{2}}f\|_{\mathcal{M}^{q,a}}\leq C \|(-\Delta)^{\frac{s}{2}}f\|_{\mathcal{M}^{p,a}}^\theta \|f\|^{1-\theta}_{\dot{B}^{-\beta, \infty}_\infty},
\end{equation}
where $1<p < q<+\infty$, $\theta=p/q$, $s_1=\theta s -(1-\theta)\beta$ and $-\beta<s_1<s$.
\end{Theoreme}

The plan of this article is the following. In Section \ref{Secc_Notation} we present our general framework which is given by stratified Lie groups. These groups are quite natural generalization of $\mathbb{R}^n$ but they present some particularities that should be taken into account in the computations.  In Section \ref{Secc_Functional_Spaces} we give the precise definition of all the functional spaces used in the previous inequalities and in Section  \ref{Secc_Proof2} we present the proof of Theorem \ref{Theorem2}. Finally, in Section \ref{Secc_Generalization} we give the proof of Theorem \ref{Theorem3} and some variations of the previous results.
\section{Stratified Lie groups: notation and basic properties}\label{Secc_Notation}
As said in the introduction, stratified Lie groups are natural generalizations of $\mathbb{R}^n$ when considering general dilation structures. Although stratified Lie groups share common features with $\mathbb{R}^n$, there are some special points that must be taken into account: for example these groups are no longer abelian and this fact requires to be carefull in some computations, furthermore from the geometric point of view, the inner geometric structure of these groups can be very different from the euclidean setting. It is then necessary to recall some basic facts about stratified Lie groups, for further information see \cite{Folland}, \cite{Folland2}, \cite{Varopoulos}, \cite{Stein2} and the references given there in. \\

We start with the notion of \textit{homogeneous group} $\mathbb{G}$ which is the data of $\mathbb{R}^{n}$ equipped with a structure of Lie group and we will always suppose that the origin is the identity. We define a \textit{dilation structure} by fixing integers $(a_{i})_{1\leq i\leq n}$ such that $1=a_{1}\leq... \leq a_{n}$ and by writing:
\begin{eqnarray}
\delta_{\alpha }:  \mathbb{R}^{n } & \longrightarrow & \mathbb{R}^{n } \label{dilat} \\
x &\longmapsto & \delta_{\alpha}[x]=(\alpha^{a_{1}}x_{1},...,\alpha^{a_{n}}x_{n}).\nonumber
\end{eqnarray}
We will often note $\alpha x$ instead of $\delta_{\alpha}[x]$ and $\alpha$ will always indicate a strictly positive real number. \\

Of course, the Euclidean space $\mathbb{R}^{n}$ with its group structure and provided with its usual dilations
(i.e. $a_{i}=1$, for $i=1,...,n$) is a homogeneous group. Here is another example: if $x=(x_{1}, x_{2}, x_{3})$ is an element of $\mathbb{R}^{3}$, we can fix a dilation by writing $\delta_{\alpha}[x]=(\alpha x_{1 }, \alpha x_{2}, \alpha^{2 } x_{3})$ for $\alpha>0$. Then, the well suited group law with respect to this dilation is given by
$ x\cdot y=(x_{1}, x_{2}, x_{3})\cdot(y_{1}, y_{2}, y_{3})=(x_{1}+y_{1}, x_{2}+y_{2}, x_{3}+y_{3}+\frac{1}{2}(x_{1}y_{2}-y_{1}x_{2 }))$. Remark in particular that this group law is no longer abelian. The triplet $(\mathbb{R}^{3}, \cdot,\delta )$ corresponds to the Heisenberg group $\mathbb{H}^{1}$ which is the first non-trivial example of a homogeneous group. The \emph{homogeneous dimension} with respect to dilation structure (\ref{dilat}) is given by $N=\displaystyle{\sum_{1\leq i\leq n}}a_{i }$. We observe that it is always larger than the topological dimension $n$ since each integer $a_{i}$ verifies $a_{i}\geq 1$ for all $i=1,...,n$. For instance, in the Heisenberg group $\mathbb{H}^{1}$ we have $N=4$ and $n=3$ while in the Euclidean case these two concepts coincide. Now we will say that a function on $\mathbb{G}\setminus \{0\}$ is \textit{homogeneous} of degree $\lambda \in \mathbb{R}$ if $f(\delta_{\alpha}[x])=\alpha^{\lambda}f(x)$ for all $\alpha>0$.  In the same way, we will say that a differential operator $D$ is homogeneous of degree $\lambda$ if $D(f(\delta_{\alpha}[x]))=\alpha^{\lambda}(Df)(\delta_{\alpha}[x])$, for all $f$ in operator's domain. In particular, if $f$ is homogeneous of degree $\lambda$ and if $D$ is a differential operator of degree $\mu$, then $Df$ is homogeneous of degree $\lambda-\mu$. The presence of a dilation structure is one of most important features of stratified Lie groups and the homogeneity with respect to these dilations will play a useful role in our computations.\\

From the point of view of measure theory, homogeneous groups behave in a traditional way since Lebesgue measure $dx$ is bi-invariant and coincides with the Haar measure, thus for any subset $E$ of $\mathbb{G}$ we will note its measure as $|E|$.
This fact also allows us to define Lebesgue spaces in a classical way (see also Section \ref{Secc_Functional_Spaces} below). The convolution will be a very useful tool in our computations, and for two functions $f$ and $g$ on $\mathbb{G}$ it is defined by
\begin{equation*}
f\ast g(x)=\int_{\mathbb{G}}f(y)g(y^{-1}\cdot x)dy=\int_{\mathbb{G}}f(x\cdot y^{-1})g(y)dy,  \quad x\in \mathbb{G}.
\end{equation*}
However, since the group law of a stratified Lie group is not necessarly commutative, we do not have in general the identity $f\ast g=g\ast f$ and we need to take care of this fact. Nevertheless, we have at our disposal Young's inequalities:
\begin{Lemme}
If $1\leq p, q, r\leq +\infty$ such that $1+\frac{1}{q}=\frac{1}{p}+\frac{1}{r}$.
If $f\in L^{p}(\mathbb{G})$ and $g\in L^{r}(\mathbb{G})$, then $f\ast g \in L^{q}(\mathbb{G})$ and we have the inequality $\|f\ast g\|_{L^q}\leq \|f\|_{L^p}\|g\|_{L^r}$.
\end{Lemme}
A proof is given in \cite{Folland2}. A weak version of Young's inequalities will be stated in Section \ref{Secc_Functional_Spaces}.\\

For a homogeneous group $\mathbb{G}=(\mathbb{R}^{n}, \cdot, \delta)$ we consider now its Lie algebra $\mathfrak{g}$ whose elements can be conceived in two different ways: as \textit{left}-invariant vector fields or as \textit{right}-invariant vector fields. The left-invariant vectors fields $(X_j)_{1\leq j\leq n}$ are determined by the formula
\begin{equation*}
(X_{j}f)(x)=\left.\frac{\partial f(x\cdot y)}{\partial y_{j}}\right|_{y=0}=\frac{\partial f}{\partial x_{j}}+\sum_{j<k}q^{k}_{j}(x)\frac{\partial f}{\partial x_{k}},
\end{equation*}
where $q^{k}_{j}(x)$ is a homogeneous polynomial of degree $a_{k}-a_{j}$ and $f$ is a smooth function on $\mathbb{G}$. By this formula one deduces easily that these vectors fields are homogeneous of degree $a_{j}$ and we have $X_{j}\left(f(\alpha x)\right)=\alpha^{a_{j}}(X_{j}f)(\alpha x)$. We will note $(Y_{j})_{1\leq j\leq n}$ the right invariant vector fields defined in a totally similar way:
$$(Y_{j}f)(x)=\left.\frac{\partial f(y\cdot x)}{\partial y_{j}}\right|_{y=0}.$$

A homogeneous group $\mathbb{G}$ is \emph{stratified} if its Lie algebra $\mathfrak{g}$ breaks up into a sum of linear subspaces
$\mathfrak{g}=\bigoplus_{1\leq j\leq k } E_{j}$ such that $E_{1}$ generates the algebra $\mathfrak{g}$ and $[E_{1}, E_{j}]=E_{j+1}$ for $1\leq j < k$ and $[E_{1}, E_{k}]=\{0\}$ and $E_{k}\neq\{0\}$, but $E_{j}=\{0\}$ if $j>k$. Here $[E_{1}, E_{j}]$ indicates the subspace of $\mathfrak{g}$ generated by the elements $[U, V]=UV-VU$ with $U\in E_{1}$ and $V\in E_{j}$. The integer $k$ is called the \emph{degree} of stratification of $\mathfrak{g}$. For example, on Heisenberg group $\mathbb{H}^1$, we have $k=2$ while in the Euclidean case $k=1$.  \\

We will suppose from now on that $\mathbb{G}$ is \textbf{stratified} with homogeneous dimension\footnote{The lower bound $N\geq 4$ corresponds to the homogeneous dimension of the Heisenberg group $\mathbb{H}^1$, which is the simplest non-trivial stratified Lie group.} $N\geq 4$. Within this framework, we will fix once and for all the family of vectors fields
\begin{equation*}
{\bf X}=\{X_1,...,X_m\},
\end{equation*}
such that $a_{1}=a_{2}=\ldots=a_{m}=1$ $(m<n)$, then the family $\textbf{X}$ is a base of $E_{1}$ and generates the Lie algebra of $\mathfrak{g}$, which is precisely the H\"ormander's condition (see \cite{Folland2} and \cite{Varopoulos}) and this particular choice ensures several important properties, in particular to the family $\textbf{X}$ is associated the Carnot-Carath\'eodory distance $d$ which is left-invariant and compatible with the topology on $\mathbb{G}$ (see \cite{Varopoulos} for more details) and for any $x\in \mathbb{G}$ we will denote by $|x|=d(x,e)$ and for $r>0$ we form open balls by writing $B(x,r)=\{y\in \mathbb{G}: d(x,y)<r\}$. By simple homogeneity arguments we obtain that stratified Lie groups have polynomial volume growth since we have $|B(\cdot, r)|=r^N |B(\cdot, 1)|$.\\

The main tools of this paper depend on the properties of the gradient, the Laplacian and the associated heat kernel, but before introducing them, we make here some remarks on general vectors fields $X_{j}$ and $Y_{j}$. Let us fix some notation: for any multi-index $I=(i_{1},...,i_{n})\in \mathbb{N}^{n}$, one defines $X^{I}$ by $X^{I}=X_{1}^{i_{1}}\dots X_{n}^{i_{n}}$ and $Y^{I}$ by $Y^{I}=Y_{1}^{i_{1}}\dots Y_{n}^{i_{n}}$, furthermore we denote by $|I|= i_{1}+\ldots+i_{n}$ the order of the derivation of the operators $X^I$ or $Y^I$ and $d(I)=a_{1}i_{1}+\ldots+a_{n}i_{n}$ the homogeneous degree of these ones. Now, for $\varphi, \psi \in \mathcal{C}^{\infty}_{0}(\mathbb{G})$ we have the equality
\begin{equation*}
\int_{\mathbb{G}}\varphi(x)(X^{I}\psi)(x)dx=(-1)^{|I|}\int_{\mathbb{G}}(X^{I}\varphi)(x)\psi(x)dx.
\end{equation*}
The interaction of operators $X^{I}$ and $Y^{I}$ with convolutions is clarified by the following identities:
\begin{equation}\label{Convolution_Left_Right}
X^{I}(f*g)=f*(X^{I}g), \qquad Y^{I}(f*g)=(Y^{I}f)*g, \qquad (X^{I}f)*g=f*(Y^{I}g).
\end{equation}
Finally, one will say that a function $f \in \mathcal{C}^{\infty}(\mathbb{G})$ belongs to the Schwartz class $\mathcal{S}(\mathbb{G})$ if the following semi-norms are bounded for all $k\in \mathbb{N}$ and any multi-index $I$:  $N_{k,I}(f)=\underset{x\in \mathbb{G}}{\sup } \, (1+|x|)^{k }|X^{I}f(x)|$.
\begin{Remarque}
\emph{To characterize the Schwartz class $\mathcal{S}(\mathbb{G})$ we can replace vector fields $X^I$ in the semi-norms $N_{k,I}$ above by right-invariant vector fields $Y^I$.}
\end{Remarque}
For a proof of these facts and for further details see \cite{Folland2} and \cite{Furioli2}.\\

We define now the \textit{gradient} on $\mathbb{G}$ from vectors fields of homogeneity degree equal to one (\textit{i.e.} those composing the family $\textbf{X}$) by fixing $\nabla = (X_{1},...,X_{m})$. This operator is of course left invariant and homogeneous of degree $1$. The length of the gradient is given by the formula $|\nabla f|= \left((X_{1}f)^{2}+... +(X_{m}f)^{2 } \right)^{1/2}$. We also define the right invariant gradient $\widetilde{\nabla}=(Y_1,...,Y_m)$, and using (\ref{Convolution_Left_Right}) we have the identity $(\nabla f)\ast g = f\ast (\widetilde{\nabla} g)$. We define now the Laplacian we are going to work with. Let us notice that in this setting there is not a single way to build a Laplacian, see for example \cite{Furioli2}. In this article we will use the Laplacian, denoted by $\mathcal{J}$, which is given from the family $\textbf{X}$ in the following way
\begin{equation}\label{Def_Laplacian}
\mathcal{J}=\nabla^{*}\nabla=-\sum_{j=1}^{m}X^{2}_{j}.
\end{equation}
This is a positive self-adjoint, hypo-elliptic operator (since the family $\textbf{X}$ satisfies the H\"ormander's condition), having as domain of definition $L^2(\mathbb{G})$. Its associated \textit{heat operator} on $\mathbb{G}\times]0, +\infty[$ is given by $\partial_{t}+\mathcal{J}$. We recall now some well-known properties of the heat operator and its associated kernel.
\begin{Theoreme}\label{Theorem_Heat_Properties} There exists a unique family of continuous linear operators $(H_{t})_{t>0}$ defined on $L^{1}+L^{\infty}(\mathbb{G})$ with the semi-group property $H_{t+s}=H_{t}H_{s}$ for all $t, s>0$ and $H_{0}=Id$, such that:
\begin{enumerate}
\item[1)] the Laplacian $\mathcal{J}$ is the infinitesimal generator of the semi-group  $H_{t}=e^{-t\mathcal{J}}$;
\item[2)] $H_{t}$ is a contraction operator on $L^{p}(\mathbb{G})$ for $1\leq p\leq +\infty$ and for $t>0$;
\item[3)] the semi-group $H_t$ admits a convolution kernel $H_{t}f=f\ast h_{t}$ where $h_{t}(x)=h(x, t) \in \mathcal{C}^{\infty}(\mathbb{G}\times]0, +\infty[)$ is the heat kernel which satisfies the following points:
\begin{enumerate}
\item $(\partial_{t}+\mathcal{J})h_{t}=0$ on $\mathbb{G}\times]0, +\infty[$, and $h(x, t)=h(x^{-1}, t)$, $h(x, t)\geq 0$ and $\displaystyle{\int_{\mathbb{G}}}h(x,t)dx=1$,
\item $h_{t}$ has the semi-group property:  $h_{t}\ast h_{s}=h_{t+s}$ for $t, s>0$ and we have $h(\delta_{\alpha}[x], \alpha^{2}t)=\alpha^{-N}h(x, t)$,
\item For every $t>0$, $x\mapsto h(x, t)$ belong to the Schwartz class in $\mathbb{G}$.
\end{enumerate}
\item[4)] For $f\in \mathcal{C}^\infty(\mathbb{G})$ and for $t>0$ we have $\mathcal{J} H_t(f)=H_t \mathcal{J} (f)$.
\end{enumerate}
\end{Theoreme}
For a detailed proof of these and other important facts concerning the heat semi-group see \cite{Folland2} and \cite{Saka}.\\

To close this section we recall the definition of the Laplacian's fractional powers. If $s>0$ we write
\begin{equation}\label{Def_Laplacian_Pos}
\mathcal{J}^s f(x)=\frac{1}{\Gamma(k-s)}\int_{0}^{+\infty}t^{k-s-1}\mathcal{J}^k H_tf(x)dt,
\end{equation}
for all $f\in \mathcal{C}^{\infty}(\mathbb{G})$ with $k$ an integer greater than $s$.  The interaction between this operator and the heat kernel is given by the following lemma.
\begin{Lemme}
If $1\leq p\leq +\infty$, for $s>0$ and for $t>0$ we have the estimate $\|\mathcal{J}^{\frac s2}h_{t}\|_{L^{p}}\leq C t^{-\frac{s+N(1-\frac 1p)}{2}}$.
\end{Lemme}
See  \cite{Saka} for a proof. 
\section{Functional spaces}\label{Secc_Functional_Spaces}

We give in this section the precise definition of all the functional spaces involved in Theorems  \ref{Theorem2} and \ref{Theorem3}. In a general way, given a norm $\|\cdot\|_{X}$, we will define the corresponding functional space $X(\mathbb{G})$ by $\{f\in \mathcal{S}'(\mathbb{G}): \|f\|_{X}<+\infty\}$. The constant that appear in this paper such as $C$ may change from one occurrence to the next.\\

\begin{enumerate}
\item[$\bullet$] \textbf{Lebesgue spaces} $L^p(\mathbb{G})$. For a measurable function $f:\mathbb{G}\longrightarrow\mathbb{R}$ and for $1\leq p <+\infty$ we define Lebesgue space by the norm $\|f\|_{L^p}=\displaystyle{\left(\int_{\mathbb{G}}|f(x)|^p dx\right)^{1/p}}$, while for $p=+\infty$ we have  $\|f\|_{L^\infty}=\underset{x\in\mathbb{G}}\esssup |f(x)|$. Let us notice that we also have the following characterization using the distribution function $\displaystyle{\|f\|^{p}_{L^{p}}=p\int_{0}^{+\infty}\alpha^{p-1 } |\{x\in \mathbb{G}:|f(x)|> \alpha\}|d\alpha}$.
\item[$\bullet$] \textbf{weak-Lebesgue spaces} $L^{p,\infty}(\mathbb{G})$. We define them  as the set of all measurable functions $f:\mathbb{G}\longrightarrow\mathbb{R}$ such that $\|f\|_ {L^{p,\infty}}=\underset{\alpha>0}{\sup}\{\alpha \cdot |\{x\in \mathbb{G}:|f(x)|> \alpha\}|^{1/p}\}$ is finite. We will need the following version of Young's inequality where weak $L^p$ spaces are involved:
\begin{Lemme}\label{Lemma_Weak_Young}
Let $p,q,r>1$. If $f\in L^{p,\infty}(\mathbb{G})$ and if $g\in L^r(\mathbb{G})$, then $f\ast g \in L^q(\mathbb{G})$ with $1+\frac{1}{q}=\frac{1}{p}+\frac{1}{r}$ and we have the inequality $\|f\ast g\|_{L^q}\leq \|f\|_{L^{p, \infty}}\|g\|_{L^r}$. 
\end{Lemme}
See a proof of this Lemma in \cite{Grafakos}, Theorem 1.4.24.
\item[$\bullet$]\textbf{Sobolev spaces} $\dot{W}^{s,p}(\mathbb{G})$. If $1<p<+\infty$ and for $s>0$ we have $\|f\|_ {\dot{W}^{s,p}} =\|\mathcal{J}^{\frac{s}{2}}f\|_{L^{p}}$, while if $p=s=1$ we will note $\|f\|_ {\dot{W}^{1,1}} =\|\nabla f\|_{L^{1}}$. We recall classical Sobolev inequalities in this setting:
\begin{equation}\label{ClassicalSobolevInequalities}
\|f\|_ {L^{\frac{N}{N-1}}} =\|\nabla f\|_{L^{1}} \qquad \mbox{and}\qquad \|f\|_ {L^{q}} =\|f\|_{\dot{W}^{s,p}}, \quad \mbox{ {\small where $1<p<q$ and $-\frac{N}{q}=s-\frac{N}{p}$}}.
\end{equation}
\item[$\bullet$]\textbf{weak Sobolev spaces} $\dot{W}^{s,p}_\infty(\mathbb{G})$. These spaces are defined just as classical Sobolev spaces, but we replace the $L^p$ norm by the weak $L^p$ one as follows:
\begin{equation*}
\|f\|_ {\dot{W}_\infty^{s, p}} =\|\mathcal{J}^{\frac s2}f\|_{L^{p,\infty}} \quad \mbox{ with } 1<p<+\infty \mbox{ and } s>0.
\end{equation*}
\item[$\bullet$] \textbf{Besov spaces} $\dot{B}^{s,q}_p(\mathbb{G})$. There are many different (and equivalent) ways to define these spaces in the setting of stratified Lie groups. In this article we will mainly use the thermic definition given by
\begin{equation*}
\|f\|_{\dot{B}^{s, q}_{p}}=\left(\int_{0}^{+\infty}t^{(m-s/2)q}\left \|\frac{\partial^{m}H_{t}f}{\partial t^{m}}(\cdot)\right\|^{q}_{L^{p}}\frac{dt}{t}\right)^{1/q},
\end{equation*}
for $1\leq p, q\leq +\infty, s>0$ and $m$ an integer such that $m>s/2$. For Besov spaces of indices $(-\beta, \infty, \infty)$ which appear in all the improved Sobolev inequalities we have:
\begin{equation*}
\|f\|_{\dot{B}^{-\beta,\infty}_{\infty}}=\underset{t>0}{\sup}\;\;t^{\beta/2 } \|H_{t}f\|_ {L^\infty}.
\end{equation*}
Recall that for $0<s<1$ we have the inequality 
\begin{equation}\label{ConstantEquivalenceBesov}
\|\mathcal{J}^{\frac s2}f\|_{\dot{B}^{-\beta-s,\infty}_{\infty}}\leq M\|f\|_{\dot{B}^{-\beta,\infty}_{\infty}},
\end{equation}
where $M$ is a universal constant. 
\item[$\bullet$]\textbf{Lorentz spaces}
 $\Lambda^p{(w)}(\mathbb{G})$. Let $f:\mathbb{G}\longrightarrow \mathbb{R}$ be a measurable function. We define $f^*$, the non-increasing rearrangement of the function $f$, by the expression $f^*(t)=\inf\{\alpha\geq 0 : |\{x\in \mathbb{G}: |f(x)|>\alpha\}|\leq t\}$. We will say that a nonnegative locally integrable function $w:\mathbb{R}^+\longrightarrow \mathbb{R}^+$ belongs to the Ari\~no-Muckenhoupt class $B_p$ for $1\leq p<+\infty$, if there exists $C>0$ such that
\[
\int_{r}^{+\infty} \left(\frac{r}{t}\right)^p w(t)dt\leq C \int_{0}^{r}w(t)dt, \mbox{ for all } 0<r<+\infty.
\]

It is not difficult to see that if $0<p<q<+\infty$, then we have the inclusion of classes $B_p\subset B_q$. We define the Lorentz spaces $\Lambda^p(w)$ with $1\leq p<+\infty$ by the formula
$$\|f\|_{\Lambda^p(w)}=\left(\int_0^{+\infty} (f^*(t))^p w(t) dt\right)^\frac{1}{p}.$$
As said in the introduction, the choice of the $B_p$ class is due to the fact that this class of weights characterizes the boundedness of the Hardy-Littlewood maximal operator $\mathcal{M}_{B}$, given for a measurable function $f$ by
\begin{equation}\label{Definition_HL_maximal_function}
\mathcal{M}_{B}f(x)=\displaystyle{\underset{B \ni x}{\sup } \;\frac{1}{|B|}\int_{B }|f(y)|dy},\quad \mbox{where } B \mbox{ is an open ball,}
\end{equation}
 on the spaces $\Lambda^p(w)$ for $1<p<+\infty$: $\|\mathcal{M}_{B}f\|_{\Lambda^p(w)}\leq C \|f\|_{\Lambda^p(w)}$, where $C$ is depending on the quantity
$$[w]_{B_p}=\sup_{r>0}\left\{r^p\, \left(\int_{r}^{+\infty}\frac{w(t)}{t^p}dt \right) \big/ \left( \int_0^r w(t)dt\right)\right\}.$$
For more properties of these weights and the associated classical Lorentz spaces see \cite{Arino}, \cite{CarroSoria}, \cite{Soria} and \cite{Carro}.
\item[$\bullet$]\textbf{Lorentz-Sobolev spaces} $\dot{\Lambda}^{s,p}(w)(\mathbb{G})$. Once we have fixed the base space $\Lambda^p(w)$, the homogeneous Lorentz-Sobolev spaces are easy to define and are given for $1<p<+\infty$ and for $s>0$ in the following way
$$\|f\|_{\dot{\Lambda}^{s,p}(w)}=\left(\int_0^{+\infty} \big( (\mathcal{J}^{\frac{s}{2}}f)^*(t)\big)^p w(t) dt\right)^\frac{1}{p}.$$

\item[$\bullet$]\textbf{weak Lorentz spaces} $\Lambda^{p,\infty}(w)(\mathbb{G})$. Let $w$ a weight in $\mathbb{R}^+$. For $0<p<+\infty$, the weak Lorentz space $\Lambda^{p,\infty}(w)$ is the class of all measurable functions $f:\mathbb{G}\longrightarrow \mathbb{R}$ such that 
$$\|f\|_{\Lambda^{p,\infty}(w)}=\sup_{t>0}f^*(t)W^{1/p}(t)<+\infty,$$
where $W(t)=\displaystyle{\int_0^t}w(s)ds$. The weak Lorentz spaces were introduced in \cite{CarroSoria} and further investigated in \cite{Carro-Soria}, \cite{Carro-Garsia del Amo-Soria} and \cite{Carro-Soria-London}. The problem of characterizing when the weak type Lorentz spaces $\Lambda^{p,\infty}(w)$, $0<p<+\infty$ are  Banach spaces was studied in \cite{Soria}.

\item[$\bullet$]\textbf{weak Lorentz-Sobolev spaces} $\dot{\Lambda}^{s,p,\infty}(w)(\mathbb{G})$. For $1<p<+\infty$ the homogeneous weak Lorentz-Sobolev spaces are given by
$$\|f\|_{\dot{\Lambda}^{s,p,\infty}(w)}=\sup_{t>0} (\mathcal{J}^{\frac{s}{2}}f)^*(t) \left(\int_0^t w(s)ds\right)^{1/p}.$$

\item[$\bullet$]\textbf{Morrey spaces} $\mathcal{M}^{p,a}(\mathbb{G})$. For $1<p<+\infty$ and $0\leq a < N$, we define Morrey spaces as the space of locally integrable functions such that
$$\|f\|_{\mathcal{M}^{p,a}}=\underset{x_0\in \mathbb{R}^n}{\sup}\; \underset{0<r<+\infty}{\sup}\left(\frac{1}{r^{a}}\int_{B(x_0,r )}|f(x)|^p dx\right)^{\frac{1}{p}}<+\infty.$$
Morrey spaces are indeed a generalization of Lebesgue spaces since when $a=0$ we have $\mathcal{M}^{p,0}\simeq L^p$. The use of Morrey and Morrey-Sobolev spaces in this article is due to the fact that the Hardy-Littlewood maximal operator (\ref{Definition_HL_maximal_function}) is also bounded in such spaces. See more details in \cite{Chiarenza} in the framework of $\mathbb{R}^n$ and \cite{Guliyev} in the setting of stratified Lie groups. See also \cite{Sawano} and the references given there in for other interesting generalizations.
\item[$\bullet$]\textbf{Morrey-Sobolev spaces} $\dot{\mathcal{M}}^{s,p,a}(\mathbb{G})$. For $0<s$ and $1<p<+\infty$ with $0\leq a < N$ we consider the homogeneous Morrey-Sobolev spaces $\dot{\mathcal{M}}^{s,p,a}$ by the quantity
$$\|f\|_{\dot{\mathcal{M}}^{s,p,a}}=\|\mathcal{J}^{\frac{s}{2}}f\|_{\mathcal{M}^{p,a}}.$$
\end{enumerate}
\section{Proof of Theorem \ref{Theorem2}}\label{Secc_Proof2}
We will prove here, in the framework of stratified Lie groups, the inequality
$$\|f\|_{\dot{\Lambda}^{s_1,q}(w)}\leq C \|f\|_{\dot{\Lambda}^{s,p}(w)}^\theta \|f\|^{1-\theta}_{\dot{B}^{-\beta, \infty}_\infty},$$
where $f:\mathbb{G}\longrightarrow\mathbb{R}$ is a function such that $f\in \dot{\Lambda}^{s,p}(w)(\mathbb{G})\cap \dot{B}^{-\beta, \infty}_\infty(\mathbb{G})$ with $1<p<q<+\infty$, $\theta=p/q$, $s_1=\theta s -(1-\theta)\beta$ and $-\beta<s_1<s$. We will always assume here that $w$ is a weight in the Ari\~no-Muckenhoupt class $B_p$. The reason for this  particular choice of weights relies on the fact that we will need the boundedness of the Hardy-Littlewood maximal operator on Lorentz $\Lambda^{p}(w)$ spaces and this is ensured by the condition $w\in B_p$. See \cite{Arino} and \cite{Carro} for details.\\

By the definition of Lorentz-Sobolev spaces given in Section \ref{Secc_Functional_Spaces}, this inequality can be rewritten in the following way
$$\|\mathcal{J}^{\frac{s_1}{2}} f\|_{\Lambda^{q}(w)}\leq C \|\mathcal{J}^{\frac{s}{2}}f\|_{\Lambda^{p}(w)}^\theta \|f\|^{1-\theta}_{\dot{B}^{-\beta, \infty}_\infty}.$$
For the proof of this inequality, we will use a variant of Hedberg's inequality. Indeed, since $0<s_1<s$, we use the characterization of the positive powers of the Laplacian given in (\ref{Def_Laplacian_Pos}) and we have for $k>s/2>s_1/2$
\begin{eqnarray*}
\mathcal{J}^{\frac{s_1}{2}}f(x)&=&\frac{1}{\Gamma(k-s_1/2)}\int_{0}^{+\infty}t^{k-\frac{s_1}{2}-1}\mathcal{J}^kH_t f(x)dt \\
&=&\frac{1}{\Gamma(k-s_1/2)}\left(\int_{0}^{T}t^{k-\frac{s_1}{2}-1}\mathcal{J}^kH_t f(x)dt+ \int_{T}^{+\infty}t^{k-\frac{s_1}{2}-1}\mathcal{J}^kH_t f(x)dt\right),
\end{eqnarray*}
where $T$ will be defined below. In particular we have
\begin{equation}\label{Hedberg_Inequality_1}
|\mathcal{J}^{\frac{s_1}{2} }f(x)|\leq \frac{1}{\Gamma(k-s_1/2)}\left(\int_{0}^{T}t^{k-\frac{s_1}{2}-1}|\mathcal{J}^kH_t f(x)|dt+ \int_{T}^{+\infty}t^{k-\frac{s_1}{2}-1}|\mathcal{J}^kH_t f(x)|dt\right).
\end{equation}
For the first integral of the right-hand side of the previous formula we will use the following fact.
\begin{Lemme}\label{Lemma_Maximal_1}
Let $f\in \mathcal{S}'(\mathbb{G})$ and $\varphi\in \mathcal{S}(\mathbb{G})$. We denote by $\mathcal{M}_{\varphi}(f)$ the maximal function of $f$ (with respect to $\varphi$) which is given by the expression
\begin{equation*}
\mathcal{M}_{\varphi}f(x)=\underset{0<t<+\infty}{\sup}\{|f\ast\varphi_{t}(x)|\}, \quad \mbox{with } \varphi_{t}(x)=t^{-N/2}\varphi(t^{-1/2}x).
\end{equation*}
If the function $\varphi$ is such that $|\varphi(x)|\leq C(1+|x|)^{-N-\varepsilon}$ for some $\varepsilon>0$, then we have the following pointwise inequality
\begin{equation*}
\mathcal{M}_{\varphi}f(x)\leq C \mathcal{M}_{B}f(x),
\end{equation*}
where $\mathcal{M}_{B}f(x)$ is the Hardy-Littlewood maximal function defined by (\ref{Definition_HL_maximal_function}).
\end{Lemme}
For a proof of this lemma see \cite{Grafakos} or \cite{Folland2}. With this lemma in mind, and since $k>s/2$, we remark that we have the identity $\mathcal{J}^kH_t f(x)=\mathcal{J}^{k-\frac{s}{2}}h_t \ast \mathcal{J}^{\frac{s}{2}} f(x)$. Now, by homogeneity we obtain $\mathcal{J}^{k-\frac{s}{2}}(h_t)(x)=t^{-k+\frac{s}{2}} \big(\mathcal{J}^{k-\frac{s}{2}}h_t\big)(x)$ and if we denote $\varphi_t$ by $\varphi_t(x)= \big(\mathcal{J}^{k-\frac{s}{2}}h_t\big)(x)$ we have that $\varphi_t(x)=t^{-N/2}\varphi(t^{-1/2}x)$, moreover, since the heat kernel $h_t$ is a smooth function, with the previous notation we obtain $|\varphi(x)|\leq C(1+|x|)^{-N-\varepsilon}$. Then we can write
$$\mathcal{J}^kH_t f(x)= t^{-k+\frac{s}{2}} \varphi_t \ast \mathcal{J}^{\frac{s}{2}} f(x),$$
and applying the Lemma \ref{Lemma_Maximal_1} we have the following pointwise inequality for the first term of (\ref{Hedberg_Inequality_1}):
\begin{equation*}
|\mathcal{J}^kH_t f(x)|= t^{-k+\frac{s}{2}}\mathcal{M}_B\left( \mathcal{J}^{\frac{s}{2}} f\right) (x).
\end{equation*}
Now, for the second integral of the right-hand side of (\ref{Hedberg_Inequality_1}) we simply use the fact that $\|\mathcal{J}^k f\|_{\dot{B}^{-\beta-2k, \infty}_{\infty}}\simeq \|f\|_{\dot{B}^{-\beta, \infty}_{\infty}}$ and the thermic definition of Besov spaces to obtain
\begin{equation*}
|\mathcal{J}^k H_tf(x)|=|H_t \mathcal{J}^k f(x)|\leq C t^{\frac{-\beta-2k}{2}}\|\mathcal{J}^k f\|_{\dot{B}^{-\beta-2k, \infty}_{\infty}}.
\end{equation*}
With these two inequalities at hand, we apply them in (\ref{Hedberg_Inequality_1}) and one has
\begin{eqnarray*}
|\mathcal{J}^{\frac{s_1}{2}}f(x)|&\leq &\frac{C}{\Gamma(k-s_1/2)}\left(\int_{0}^{T} t^{k-\frac{s_1}{2}-1} t^{-k+\frac{s}{2}}\mathcal{M}_B\left( \mathcal{J}^{\frac{s}{2}} f\right) (x)dt+ \int_{T}^{+\infty}t^{k-\frac{s_1}{2}-1} t^{\frac{-\beta-2k}{2}}\|\mathcal{J}^k f\|_{\dot{B}^{-\beta-2k, \infty}_{\infty}}dt\right)\\
&\leq & \frac{C}{\Gamma(k-s_1/2)}\left( T^{\frac{s-s_1}{2}} \mathcal{M}_B \left(\mathcal{J}^{\frac{s}{2}} f\right)(x)+ T^{\frac{-\beta-s_1}{2}}\|\mathcal{J}^k f\|_{\dot{B}^{-\beta-2k, \infty}_{\infty}}\right).
\end{eqnarray*}
We fix now the parameter $T$ by the condition
$$T=\left(\frac{\|\mathcal{J}^k f\|_{\dot{B}^{-\beta-2k, \infty}_{\infty}}}{\mathcal{M}_B \left(\mathcal{J}^{\frac{s}{2}} f\right)(x)}\right)^{ \frac{2}{\beta+s}},$$
and we obtain the following inequality
$$|\mathcal{J}^{\frac{s_1}{2}}f(x)|\leq \frac{C}{\Gamma(k-s_1/2)}\mathcal{M}_B \left(\mathcal{J}^{\frac{s}{2}} f\right)^{1-\frac{s-s_1}{\beta+s}}(x)\|\mathcal{J}^kf\|^{\frac{s-s_1}{\beta+s}}_{\dot{B}^{-\beta-2k, \infty}_{\infty}}.$$
Since $\frac{s-s_1}{\beta+s}=1-\theta$ and using again the fact $\|\mathcal{J}^k f\|_{\dot{B}^{-\beta-2k, \infty}_{\infty}}\simeq \|f\|_{\dot{B}^{-\beta, \infty}_{\infty}}$ we have
\begin{equation}\label{Pointwise_Inequality}
|\mathcal{J}^{\frac{s_1}{2}}f(x)|\leq \frac{C}{\Gamma(k-s_1/2)}\mathcal{M}_B \left(\mathcal{J}^{\frac{s}{2}} f\right)^{\theta}(x)\|f\|^{1-\theta}_{\dot{B}^{-\beta, \infty}_{\infty}}.
\end{equation}
Once we have obtained this pointwise inequality, we will use the following properties of the non-increasing rearrangement function.
\begin{Lemme} If $f,g:\mathbb{G}\longrightarrow \mathbb{R}$ are two measurable functions, we have
\begin{itemize}
\item[1)] if $|g|\leq |f|$ a.e. then $g^*\leq f^*$,
\item[2)] if $0<\theta$, then $(|f|^\theta)^*=(f^*)^\theta$.
\end{itemize}
\end{Lemme}
For a proof see Proposition 1.4.5 of \cite{Grafakos}. Recalling that $\theta=p/q$ and applying these facts to the inequality (\ref{Pointwise_Inequality}) we obtain
\begin{equation}\label{pointwise ineq}
\left((\mathcal{J}^\frac{s_1}{2}f)^*(t)\right)^q\leq C\left((\mathcal{M}_B \left(\mathcal{J}^{\frac{s}{2}} f\right))^*(t)\right)^p \|f\|^{q-p}_{\dot{B}_{\infty}^{-\beta,\infty}}.
\end{equation}
Multiplying the previous inequality by a weight $w$ from the Ari\~no-Muckenhoupt class $B_p$ and integrating with respect to the variable $t$ we obtain
\begin{equation*}
\int_{0}^{+\infty}\left((\mathcal{J}^\frac{s_1}{2}f)^*(t)\right)^q w(t)dt \leq C\int_{0}^{+\infty}\left((\mathcal{M}_B \left(\mathcal{J}^{\frac{s}{2}} f\right))^*(t)\right)^p w(t)dt\; \|f\|^{q-p}_{\dot{B}_{\infty}^{-\beta,\infty}},
\end{equation*}
and then, by the definition of classical Lorentz spaces given in Section \ref{Secc_Functional_Spaces} we have
$$\|\mathcal{J}^{\frac{s_1}{2}} f\|_{\Lambda^q(w)}\leq C\|\mathcal{M}_B(\mathcal{J}^{\frac{s}{2}} f)\|^\theta _{\Lambda^p(w)}\|f\|^{1-\theta}_{\dot{B}_{\infty}^{-\beta,\infty}}.$$
Now, since the weight $w$ belongs to the class $B_p$ with $1<p<+\infty$, we have that the Hardy-Littlewood maximal operator is bounded on the space $\Lambda^p(w)$ and we obtain
$$\|\mathcal{M}_B(\mathcal{J}^{\frac{s}{2}} f)\|_{\Lambda^p(w)}\leq \|\mathcal{J}^{\frac{s}{2}} f\|_{\Lambda^p(w)},$$
and finally we have the desired inequality for classical Lorentz spaces:
$$\|\mathcal{J}^{\frac{s_1}{2}} f\|_{\Lambda^q(w)}\leq C\|\mathcal{J}^{\frac{s}{2}}f\|_{\Lambda^p(w)}^\theta \|f\|^{1-\theta}_{\dot{B}_{\infty}^{-\beta,\infty}}.$$
\hfill $\blacksquare$\\

Now we will state in the following corollaries some interesting consequences of this previous theorem.
\begin{Corollaire}\label{Corollary_Weak_Lorentz}
Let $w\in B_p$ be a weight and let $f:\mathbb{G}\longrightarrow\mathbb{R}$ be a function such that $f\in \dot{\Lambda}^{s,p,\infty}(w)(\mathbb{G})\cap \dot{B}^{-\beta, \infty}_\infty(\mathbb{G})$. Then we have the following version of improved Sobolev inequalities of weak type:
$$ \|f\|_{\dot{\Lambda}^{s_1,q,\infty}(w)}\leq C\|f\|_{\dot{\Lambda}^{s,p,\infty}(w)}^\theta \|f\|^{1-\theta}_{\dot{B}_{\infty}^{-\beta,\infty}},$$
where $1<p < q<+\infty$, $\theta=p/q$, $s_1=\theta s -(1-\theta)\beta$ and $-\beta<s_1<s$.
\end{Corollaire}
\textbf{\textit{Proof}}. We start again with the pointwise inequality \eqref{pointwise ineq}:
$$\left((\mathcal{J}^\frac{s_1}{2}f)^*(t)\right)^q\leq C\left((\mathcal{M}_B \left(\mathcal{J}^{\frac{s}{2}} f\right))^*(t)\right)^p \|f\|^{q-p}_{\dot{B}_{\infty}^{-\beta,\infty}}.$$
Now, we multiply both parts of this inequality by $W(t)$ and we take the supremum in the variable $t$:
\begin{eqnarray*}
\|\mathcal{J}^{\frac{s_1}{2}} f\|^q_{\Lambda^{q,\infty}(w)}=\sup_{t>0}W(t)\left((\mathcal{J}^{\frac{s_1}{2}} f)^*(t)\right)^q &\leq &C \sup_{t>0}\left\{\left(\mathcal{M}_B(\mathcal{J}^{\frac{s}{2}} f)^*(t)\right)^pW(t)\right\}\|f\|^{q-p}_{\dot{B}_{\infty}^{-\beta,\infty}}\\
&\leq & C\|\mathcal{M}_B(\mathcal{J}^{\frac{s}{2}} f)\|^p_{\Lambda^{p,\infty}(w)} \|f\|^{q-p}_{\dot{B}_{\infty}^{-\beta,\infty}},
\end{eqnarray*}
since it is known (see e.g. \cite{Soria}) that for $w\in B_p$ the Hardy-Littlewood maximal operator $\mathcal{M}_B$ is bounded on $\Lambda^{p,\infty}(w)$, therefore we obtain that
$$\|\mathcal{J}^{\frac{s_1}{2}} f\|_{\Lambda^{q,\infty}(w)}\leq C\|\mathcal{J}^{\frac{s}{2}}f\|_{\Lambda^{p,\infty}(w)}^\theta \|f\|^{1-\theta}_{\dot{B}_{\infty}^{-\beta,\infty}}.$$
\hfill $\blacksquare$\\

Now we will study other variations of the previous results by considering a different type of weights. To be more precise, we will study two-weighted inequalities and in what follows, for $v$ and $w$ two weigths and for $t>0$, we will denote by  $V(t)$ and $W(t)$ the quantities $V(t)=\displaystyle{\int_0^t} v(s) ds$ and $W(t)=\displaystyle{\int_0^t} w(s) ds$.\\

Our first two-weighted improved Lorentz-Sobolev inequality is given in the following corollary.
\begin{Corollaire} Let $1<p<q<+\infty$ and let $(v,w)$ be a pair of positive weights satisfying the following properties
$$\sup_{t>0} \frac{W(t)^{1/p}}{V(t)^{1/p}}<+\infty\quad \text{ and } \quad \sup_{t>0}\left(\int_t^{+\infty} \frac{w(s)}{s^p}ds\right)^{1/p}\left(\int_0^t \frac{v(s)s^{p'}}{V(s)^{p'}}ds\right)^{1/p'}<+\infty.$$
If $f:\mathbb{G}\longrightarrow \mathbb{R}$ is a function such that $f\in \dot{\Lambda}^{s,p}(v)\cap \dot{B}_{\infty}^{-\beta,\infty} $ with $s>0$, then we have a two-weighted version of improved Sobolev inequalities
 $$ \| f\|_{\dot{\Lambda}^{s_1,q}(w)}\leq C\|f\|_{\dot{\Lambda}^{s,p}(v)}^\theta \|f\|^{1-\theta}_{\dot{B}_{\infty}^{-\beta,\infty}},$$
where $1<p < q<+\infty$, $\theta=p/q$, $s_1=\theta s -(1-\theta)\beta$ and $-\beta<s_1<s$.
\end{Corollaire}
This inequality is interesting since it is possible, under some hypotheses, to consider different weights in the left-hand side and in the right-hand side of the inequality.\\[2mm]
\textbf{\textit{Proof}}. Using the pointwise inequality \eqref{pointwise ineq} and the fact that the Hardy-Littlewood maximal operator 
$$\mathcal{M}_B:\Lambda^p(v) \longrightarrow \Lambda^p(w)$$ 
is bounded for such weights (see \cite{Sawyer} for details) we obtain the desired inequality. \hfill $\blacksquare$\\

If we are allowed to change the weights that define the Lorentz spaces in the previous inequalities, it is then also possible to change, with specific conditions on the weights, the parameters of these spaces. In the following corollary we gather some results where we consider different Lorentz spaces in the right-hand side of the inequality. Indeed, the first point is a generalization of the previous corollary and we will consider in the right-hand side Lorentz-Sobolev spaces of type $\dot{\Lambda}^{s,q_0}(v)$ instead of $\dot{\Lambda}^{s,p}(v)$ where $1<q_{0}\leq p<+\infty$. The second point allows us to study the case when $1<p<q_{0}<+\infty$ and finally, the third point treats the case when  $0<q_{0}<1$.
\begin{Corollaire} Let $0<q_{0}<+\infty$,  $s>0$, let $f:\mathbb{G}\longrightarrow\mathbb{R}$ be a measurable function and let $(v,w)$ be a pair of weights.
\begin{enumerate}
\item[1)] If $1<q_0\le p<+\infty$ and if $(v,w)$ are satisfying the following conditions 
\begin{equation}\label{two weights 1}
\sup_{t>0} \frac{W(t)^{1/p}}{V(t)^{1/{q_0}}}<+\infty
\end{equation}
and
\begin{equation}\label{two weights 2}
\sup_{t>0} \left(\int_0^t \frac{w(s)}{s^p} ds\right)^{1/p}\left(\int_0^t \frac{v(s)s^{q_0'}}{V(s)^{q_0'}}\right)<+\infty,
\end{equation}
then, if $f\in \dot{\Lambda}^{s,q_0}(v)(\mathbb{G})\cap \dot{B}^{-\beta, \infty}_\infty(\mathbb{G})$, we have the following inequality
$$\| f\|_{\dot{\Lambda}^{s_1,q}(w)}\leq C\|f\|_{\dot{\Lambda}^{s,q_0}(v)}^\theta \|f\|^{1-\theta}_{\dot{B}_{\infty}^{-\beta,\infty}},$$
where $1<p < q<+\infty$, $\theta=p/q$, $s_1=\theta s -(1-\theta)\beta$ and $-\beta<s_1<s$.

\item[2)] If $1<p<q_0<+\infty$ and $(v,w)$ are satisfying 
\begin{equation*}
\left(\int_0^{+\infty}\left(\frac{W(s)}{V(s)}\right)^{r/{q_0}}w(s)ds\right)^{1/r}<+\infty
\end{equation*}
and
\begin{equation*}
\left(\int_0^{+\infty}\left[\left(\int_s^{+\infty}\frac{w(t)}{t^p}dt\right)^{1/p}\left(\int_0^t\frac{v(t)t^{{q_0}'}}{V(t)^{{q_0}'}}dt\right)^{1/{p'}}\right]^r\frac{v(s)s^{{q_0}'}}{V(s)^{{q_0}'}}ds\right)^{1/r}<+\infty,
\end{equation*}
where $r$ is given by $\frac{1}{r}=\frac{1}{p}-\frac{1}{q}$ and $\frac{1}{q_0}+\frac{1}{{q_0}'}=1$. Then, if $f\in \dot{\Lambda}^{s,q_0}(v)(\mathbb{G})\cap \dot{B}^{-\beta, \infty}_\infty(\mathbb{G})$, we have   
$$\| f\|_{\dot{\Lambda}^{s_1,q}(w)}\leq C\|f\|_{\dot{\Lambda}^{s,q_0}(v)}^\theta \|f\|^{1-\theta}_{\dot{B}_{\infty}^{-\beta,\infty}},$$
where $1<p < q<+\infty$, $\theta=p/q$, $s_1=\theta s -(1-\theta)\beta$ and $-\beta<s_1<s$.

\item[3)] If $0<q_0<1$ and $1< p<+\infty$ and if $(v,w)$ are satisfying \eqref{two weights 1} and
\begin{equation*}
	\sup_{t>0} \frac{t}{V(t)^{1/{q_0}}}\left(\int_t^{+\infty} \frac{w(s)}{s^p} ds\right)^{1/p}<+\infty ,
\end{equation*} 
then, assuming that $f\in \dot{\Lambda}^{s,q_0}(v)(\mathbb{G})\cap \dot{B}^{-\beta, \infty}_\infty(\mathbb{G})$, we obtain
$$\| f\|_{\dot{\Lambda}^{s_1,q}(w)}\leq C\|f\|_{\dot{\Lambda}^{s,q_0}(v)}^\theta \|f\|^{1-\theta}_{\dot{B}_{\infty}^{-\beta,\infty}},$$
where $1<p < q<+\infty$, $\theta=p/q$, $s_1=\theta s -(1-\theta)\beta$ and $-\beta<s_1<s$.
\end{enumerate}		
\end{Corollaire}
\textbf{\textit{Proof}.}  From the pointwise inequality \eqref{pointwise ineq} we obtain that
$$\|\mathcal{J}^{\frac{s_1}{2}} f\|_{\Lambda^{q}(w)}\leq C\|\mathcal{M}_{B}(\mathcal{J}^{\frac{s}{2}}f)\|_{\Lambda^{p}(w)}^\theta \|f\|^{1-\theta}_{\dot{B}_{\infty}^{-\beta,\infty}}.$$
Now, under all these hypotheses on the weights $v$ and $w$, we have that the Hardy-Littlewood maximal operator $\mathcal{M}_B: \Lambda^{q_{0}}(v)\longrightarrow\Lambda^p(w)$ is bounded (see \cite{Sawyer} and \cite{Carro-Soria}) and then we obtain
$$ \| f\|_{\dot{\Lambda}^{s_1,q}(w)}\leq C\|f\|_{\dot{\Lambda}^{s,q_{0}}(v)}^\theta \|f\|^{1-\theta}_{\dot{B}_{\infty}^{-\beta,\infty}}.$$
\hfill $\blacksquare$\\

We have also the following two-weighted version of improved Sobolev inequalities of weak type:
\begin{Corollaire}
Let $1<p<+\infty$, $0<q_{0}<+\infty$. Let $(v,w)$ be a pair of weights such that
\begin{equation}\label{condweak}
\sup_{t>0}\frac{W(t)^{1/p}}{t}\int_0^tV^{-1/{q_0}}(s)ds<+\infty,
\end{equation}
and let $f:\mathbb{G}\longrightarrow\mathbb{R}$ be a function such that $f\in \dot{\Lambda}^{s, q_0,\infty}(v)(\mathbb{G})\cap \dot{B}^{-\beta, \infty}_\infty(\mathbb{G})$. Then we have the following inequality
$$\|f\|_{\dot{\Lambda}^{s_1,q,\infty}(w)}\leq C\|f\|_{\dot{\Lambda}^{s,q_0,\infty}(v)}^\theta \|f\|^{1-\theta}_{\dot{B}_{\infty}^{-\beta,\infty}},$$
where $0<q_0<+\infty$, $1<p < q<+\infty$, $\theta=p/q$, $s_1=\theta s -(1-\theta)\beta$ and $-\beta<s_1<s$.
\end{Corollaire}
\textbf{\textit{Proof}}. It is enough to follow the same lines of the Corollary \ref{Corollary_Weak_Lorentz} to obtain
$$\|\mathcal{J}^{\frac{s_1}{2}} f\|^q_{\Lambda^{q,\infty}(w)}\leq C\|\mathcal{M}_B(\mathcal{J}^{\frac{s}{2}} f)\|^p_{\Lambda^{p,\infty}(w)} \|f\|^{q-p}_{\dot{B}_{\infty}^{-\beta,\infty}},$$
since the pair of weights $(v,w)$ satisfies the condition \eqref{condweak} it implies that the operator $\mathcal{M}_B:\Lambda^{q_{0},\infty}(v) \longrightarrow \Lambda^{p,\infty}(w)$ is bounded (see \cite{Soria}) and we obtain
$$\|\mathcal{J}^{\frac{s_1}{2}} f\|^q_{\Lambda^{q,\infty}(w)} \leq C\|\mathcal{J}^{\frac{s}{2}} f\|^p_{\Lambda^{q,\infty}(v)} \|f\|^{q-p}_{\dot{B}_{\infty}^{-\beta,\infty}},$$
which is the desired inequality. \hfill $\blacksquare$
\section{Generalizations}\label{Secc_Generalization}
In this section we give some generalizations of Theorems \ref{Theorem2} and we prove Theorem \ref{Theorem3}. These generalizations are made possible since the techniques developed in our proofs are based on general harmonic analysis arguments and since many of the tools used in this article are available in other frameworks. Indeed, the spectral theory associated to the Laplace operator, the boundedness of the Hardy-Littlewood maximal operator and the use of appropiate weights in order to define well suited functional spaces are intensively studied and many interesting properties were generalized to different settings. 

\subsection{Morrey spaces}
We prove now Theorem \ref{Theorem3} in the setting of stratified Lie groups. Morrey spaces were studied in this framework by many authors, see for example the articles \cite{Guliyev}, \cite{Nakai1} and the references there in.\\

As said in the introduction, once we have at our disposal the fact that the Hardy-Littlewood maximal operator is bounded in the convenient functional framework, it is possible to improve Sobolev inequalities in the following way. The starting point of our proof is the pointwise inequality (\ref{Pointwise_Inequality}): 
$$|\mathcal{J}^{\frac{s_1}{2}}f(x)|\leq \frac{C}{\Gamma(k-s_1/2)}\mathcal{M}_B \left(\mathcal{J}^{\frac{s}{2}} f\right)^{\theta}(x)\|f\|^{1-\theta}_{\dot{B}^{-\beta, \infty}_{\infty}}.
$$
Since $\theta=p/q$ we have for $r>0$ and for $0\leq a <N$ the inequalities
\begin{eqnarray*}
\frac{1}{r^a}\int_{B(x_0, r)}|\mathcal{J}^{\frac{s_1}{2}}f(x)|^q dx &\leq & C \left(\frac{1}{r^a}\int_{B(x_0, r)}\mathcal{M}_B \left(\mathcal{J}^{\frac{s}{2}} f\right)^{p}(x)dx \right)\|f\|^{q(1-\theta)}_{\dot{B}^{-\beta, \infty}_{\infty}}\\
\left(\frac{1}{r^a}\int_{B(x_0, r)}|\mathcal{J}^{\frac{s_1}{2}}f(x)|^q dx\right)^{1/q} &\leq & C \left(\frac{1}{r^a}\int_{B(x_0, r)}\mathcal{M}_B \left(\mathcal{J}^{\frac{s}{2}} f\right)^{p}(x)dx \right)^{1/q}\|f\|^{1-\theta}_{\dot{B}^{-\beta, \infty}_{\infty}},
\end{eqnarray*}
from which we derive the estimate
$$\|\mathcal{J}^{\frac{s_1}{2}}f\|_{\mathcal{M}^{q,a}}\leq C\|\mathcal{M}_B \left(\mathcal{J}^{\frac{s}{2}} f\right)\|_{\mathcal{M}^{p,a}}^{\theta}\|f\|^{1-\theta}_{\dot{B}^{-\beta, \infty}_{\infty}}.$$
In order to conclude, we use the fact that the Hardy-Littlewood maximal operator is bounded in Morrey spaces and we obtain
$$\|\mathcal{J}^{\frac{s_1}{2}}f\|_{\mathcal{M}^{q,a}}\leq C\|\mathcal{J}^{\frac{s}{2}}f\|_{\mathcal{M}^{p,a}}^{\theta}\|f\|^{1-\theta}_{\dot{B}^{-\beta, \infty}_{\infty}},$$
which is the desired inequality stated in Theorem \ref{Theorem3}. 

\begin{Remarque}
The boundedness of the Hardy-Littlewood maximal operator was studied for generalized Morrey spaces in \cite{Guliyev}, \cite{Nakai} and \cite{Sawano}. As long as this boundedness property is satisfied it should be possible to generalize Theorem \ref{Theorem3}. Indeed, from the pointwise inequality (\ref{Pointwise_Inequality}) it should be easy (taking into account the necessary precautions) to reconstruct the corresponding norms in order to obtain an improved Sobolev-like inequality. 
\end{Remarque}
\subsection{Nilpotent Lie groups}
We consider now a more general framework than the one given by stratified Lie groups. Indeed, going one step further in the process of generalization, it is possible to consider nilpotent Lie groups since all the tools used in the proof of Theorem \ref{Theorem2} are available in these settings. \\

We recall for the sake of completness this framework. Let $\mathbb{G}$ be a connected unimodular Lie group endowed with its Haar measure $dx$. Denote by $\mathfrak{g}$ the Lie algebra of $\mathbb{G}$ and consider a family (that will be fixed from now on) of left-invariant vector fields on $\mathbb{G}$
\begin{equation*}
{\bf X}=\{X_1,...,X_k\},
\end{equation*}
satisfying the \textit{H\"ormander condition}\footnote{which means that the Lie algebra generated by the family $\textbf{X}$ is $\mathfrak{g}$.}. We endow the group $\mathbb{G}$ with a metric structure by considering the Carnot-Carath\'eodory metric associated with ${\bf X}$. See \cite{Varopoulos} for details. We will denote $\|x\|$ the distance between the origin $e$ and $x$ and $\|y^{-1}\cdot x\|$ the distance between $x$ and $y$. For $r>0$ and $x\in \mathbb{G}$, denote by $B(x,r)$ the open ball with respect to the Carnot-Carath\'eodory metric centered in $x$ and of radius $r$, and by $V(r)=\displaystyle{\int_{B(x,r)}dx}$ the Haar measure of any ball of radius $r$. When $0<r<1$, there exists $d\in \mathbb{N}^\ast$, $c_l$ and $C_l>0$ such that, for all $0<r<1$ we have
\begin{equation*}
c_l r^d\leq V(r)\leq C_l r^d.
\end{equation*}
The integer $d$ is the \textit{local} dimension of $(\mathbb{G}, {\bf X})$. When $r\geq 1$, two situations may occur, independently of the choice of the family ${\bf X}$: either $\mathbb{G}$ has polynomial volume growth and there exist $D\in \mathbb{N}^\ast$,  $c_\infty$ and $C_\infty>0$ such that, for all $r\geq 1$ we have
\begin{equation*}
c_\infty r^D\leq V(r)\leq C_\infty r^D,
\end{equation*}
or $\mathbb{G}$ has exponential volume growth, which means that there exist  $c_e, C_e, \alpha, \beta>0$ such that, for all $r\geq 1$ we have
$$c_e e^{\alpha r}\leq V(r)\leq C_e  e^{\beta r}.$$
When $\mathbb{G}$ has polynomial volume growth, the integer $D$ is called the \textit{dimension at infinity} of $\mathbb{G}$. Recall that nilpotent groups have polynomial volume growth and that a \textit{strict} subclass of the nilpotent groups consists of stratified Lie groups where $d=D$.\\

Once we have fixed the family $\textbf{X}$, we define the gradient on $\mathbb{G}$ by $\nabla = (X_1,...,X_k)$ and we consider a Laplacian $\mathcal{J}$ on $\mathbb{G}$ defined in the same way as in (\ref{Def_Laplacian})
\begin{equation*}
\mathcal{J}=-\sum_{j=1}^k X_j^2,
\end{equation*}
which is a positive self-adjoint, hypo-elliptic operator since ${\bf X}$ satisfies the H\"ormander's condition, see \cite{Varopoulos}. Its associated heat operator on $]0, +\infty[\times\mathbb{G}$ is given by $\partial_{t}+\mathcal{J}$ and we will denote by $(H_t)_{t> 0}$ the semi-group obtained from the Laplacian $\mathcal{J}$. It is worth noting that many of the properties given in Theorem \ref{Theorem_Heat_Properties} remain true for the heat semi-group $H_t$ in this general setting. For more details concerning nilpotent Lie groups see the books \cite{Varopoulos}, \cite{Folland2}, \cite{Stein2} and the articles \cite{Furioli2}, \cite{Saka}, \cite{Chamorro2} and the references there in. Fractional powers of the Laplacian can be defined in a completely similar way using the expression (\ref{Def_Laplacian_Pos}). It is then possible to define all the functional spaces given in Section \ref{Secc_Functional_Spaces} in the framework of nilpotent Lie groups. \\

With all these preliminaries, we see that we have at our disposal all the ingredients needed in order to perform the computations done in Sections  \ref{Secc_Proof2}, and thus Theorem \ref{Theorem2} can be generalized to the setting of nilpotent Lie groups.\\

\textbf{Acknowledgments.} A part of this work was performed while the second and third authors visited the University of Evry Val d'Essonne. We express our gratitude to the Laboratoire de Math\'ematiques et Mod\'elisation d'Evry (LaMME) of the University of Evry Val d'Essonne for the hospitality and excellent conditions. The second named author was partially supported by POSDRU/159/1.5/S/137750.

{\footnotesize
\begin{tabular}{lll }
Diego \textsc{Chamorro}						& \hspace{0.5cm} & Anca \textsc{Marcoci} \& Liviu \textsc{Marcoci}\\[3mm]
Laboratoire de Math\'ematiques et Mod\'elisation d'Evry	& & Department of Mathematics and Computer Science\\
(LaMME) UMR 8071, Universit\'e d'Evry Val d'Essonne		& & Technical University of Civil Engineering\\
23 Boulevard de France, 91037 Evry Cedex					& & Bucharest, Bld. Lacul Tei, no. 124, sector 2.\\
 France 			 										& & Romania\\
diego.chamorro@univ-evry.fr								& & anca.marcoci@utcb.ro \quad  lmarcoci@instal.utcb.ro
\end{tabular}
}
\end{document}